\documentclass[preprint,review,12pt]{elsarticle}
\usepackage{amsmath,amssymb,amsthm}
\usepackage{amsfonts}
\usepackage{mathrsfs}
\usepackage{CJK}

\newtheorem{Th}{Theorem}[section]
 \newtheorem{Lem}{Lemma}[section]
\newtheorem{Cor}{Corollary}[section]

\newtheorem{Def}{Definition}[section]
\newtheorem{Rem}{Remark}[section]
\newtheorem{Exa}{Example}[section]
\newtheorem{Prob}{Problem}[section]
\newtheorem{Method}{Method}[section]

\def\R{{\bf R}}

 \def\C{{\bf C}}
 
 \def\ux{\underline{x}}
 \def\ur{\underline{r}}
 \def\uv{\underline{v}}
 \def\uy{\underline{y}}

\def\i{{\bf i}}
\def\e{{\bf e}}

 \def\C{{\bf C}}

 \def\ux{\underline{x}}
 
 \def\un{\underline{n}}

\def\uxi{\underline{\xi}}
 \def\ut{\underline{t}}

\journal{Multidimensional Systems and Signal Processing}

\date{}

\begin{document}
\begin{frontmatter}

\title{\bf Edge Detection Methods Based on Modified Differential Phase Congruency of Monogenic Signal}

\author{Yan Yang}
\address{School of Mathematics (Zhuhai), Sun Yat-Sen University, \\
Zhuhai, 519080, China.\\ Email:  mathyy@sina.com}

\author{Kit Ian Kou}
\address{Department of Mathematics, University of Macau, Macao (Via Hong
Kong).\\ Email: kikou@umac.mo}

\author{Cuiming~Zou}
\address{Department of Mathematics, University of Macau, Macao (Via Hong
Kong).\\ Email: zoucuiming2006@163.com}

\begin{abstract}
Monogenic signal is regarded as a generalization of analytic
signal from one dimensional to higher dimensional space, which has been received considerable attention in the literature. It is defined by an original signal with its isotropic Hilbert transform (the combination of Riesz transform). Like the analytic signal, monogenic signal can be written in the polar form. Then it provides the signal features representation, such as the local attenuation and the local phase vector.
The aim of the paper is twofold: first, to analyze the relationship between the local phase vector and the local attenuation in the higher dimensional spaces. Secondly, a study on image edge detection using modified differential phase congruency is presented. Comparisons with competing methods on real-world images consistently show the superiority of the proposed methods.
\end{abstract}


\end{frontmatter}

\section{Introduction}\label{S1}

In the scale-space literature, there are a lot of
papers discussing Gaussian scale-space as the only linear
scale-space \cite{BWBD, I, L}. The Gaussian scale space is obtained as the solution of the heat
equation. In \cite{FS2}, M. Felsberg
and G. Sommer proposed a new linear scale-space which is generated by
the Poisson kernel, it is the so-called Poisson scale-space in
two-dimensional (2D) spaces.

The Poisson scale-space is obtained by Poisson filtering (the
convolution of the original signal and the Poisson kernel). The
harmonic conjugate (the convolution of the original signal and the
conjugate Poisson kernel) yields the corresponding figure flow at
all scales. The Poisson scale-space and its corresponding figure
flow form the monogenic scale-space \cite{FS2}. In mathematics, monogenic
scale-space is the Hardy space in the upper half complex plane. The
boundary value of a monogenic function in the upper half space is the
monogenic signal. The monogenic scale signal gives deeper insight to low level image processing.

Monogenic signal is regarded as a generalization of analytic
signal from one dimensional space to higher dimensional case, which is first studied by M.
Felsberg and G. Sommer in 2001 \cite{FS1}.
It is defined by an original signal with its Riesz transform in higher dimensions.
Under certain assumptions, monogenic function can be representation in the polar form
 and then it provides the signal features, such as the local attenuation and the
local phase vector \cite{FS1, YQS, YDQ, F}. In \cite{YQS}, we first defined the scalar-valued phase derivative
(local frequency) of a multivariate signal in higher dimensions.
Then we studied the applications in signal processing \cite{YDQ, YK, LZa, LZb, LZc}. 

Monogenic signals at any scale $s>0$ form monogenic scale-space.
The representation of monogenic scale-space is just a monogenic function in the upper half-space.
Therefore, considering the monogenic scale space with scale $s$ instead of monogenic signals, provides us more analysis tools.
 In the monogenic scale-space,
the important features in image processing, such as local phase-vector,
and local attenuation (the log of local amplitude) involving through scale $s$ are given in \cite{FS2}.
The relationship between the local attenuation
and the local phase in the intrinsically 1D cases are derived in \cite{FS2}. However, the problem is open if the signal is
not intrinsically 1D signal.

The contributions of this paper are summarized as follows.

\begin{itemize}
\item[1.] We give the solution of the problem: if the higher dimensional signal is not intrinsically 1D signal, we derive the relationship between the local phase-vector and and local attenuation.

\item[2.] We proposed the local attenuation (LA) method for edge detection operator. We establish the theoretical and experiment results on the newly methods.

\item[3.] We proposed the modified differential phase congruency (MDPC) method for the edge detection operator. We establish the theoretical and experiment results on the newly methods.

\item[4.] We show that in higher dimensional space, the instantaneous frequency in higher dimensional spaces defined by
 is equal to the minus of the scale derivative of the
local attenuation.

\item[5.] We show that the zero points of the differential phase congruency is {\bf not} equal to the extrema of the local attenuation. The nonzero extra term $$-{\rm
Vec} \left[ \left(\underline{D}\frac{\underline{v}}{|\underline{v}|} \right)\frac{\underline{v}}{|\underline{v}|} \right]\sin^2
\theta+(\sin \theta \cos \theta - \theta){\partial \over \partial s} \left({\uv \over |\ux|} \right)$$ appears in high dimensional cases.

\end{itemize}

The rest of this paper is organized as follows. In order to make it self-contained, Section \ref{S2} gives a brief introduction to some general definitions and basic properties of Hardy space, analytic signal, Clifford algebra, monogenic signal and monogenic scale space. In Section \ref{S3} we derive the relationship between the local phase-vector and and local attenuation. Various edge detection methods are provided in Section \ref{S4}. Finally, experiment results are drawn in Section \ref{S5}.

\section{Preliminaries}\label{S2}
In the present section, we begin by reviewing some definitions and basic properties of analytic signal and Hardy space \cite{Co, G1981, H1996}.

\subsection{Analytic Signal and Hardy Space}
\begin{Def}[Analytic Signal] For a square integrable real-valued function $f$, the complex-valued signal $f_A$ whose imaginary part is the Hilbert transform of its real part is called the {\it analytic signal}. That is, $$f_A(x) :=f(x)+\i \mathcal{H}(f)(x),$$
where $\mathcal{H}(f)(x)$ is the Hilbert transform (HT) of $f$ defined by
\begin{eqnarray*}
\mathcal{H}[f](x) := \frac{1}{\pi}{\rm p.v.}\int_{-\infty}^{+\infty}\frac{f(s)}{x-s}ds=\frac{1}{\pi} \lim\limits_{\varepsilon\to 0^+} \int_{\varepsilon \leq |x-s|} \frac{f(s)}{x-s}ds,\label{liu11}
\end{eqnarray*} provided this integral exists as a principal value (${\rm p.v.}$ means the Cauchy principle value).
\end{Def}

Due to its definition, the real $u$ and imaginary parts $v$ of analytic signal $f_A=u+\i v$ form the {\it Hilbert transform pairs}
\begin{eqnarray}\label{HTP}
\mathcal{H}[u]=v.
\end{eqnarray}

To proceed the properties of analytic signal, we introduce the notion of Hardy space \cite{G1981, H1996}, we will notice that the class of analytic signals is the class of boundary values of Hardy space functions.

\begin{Def}[Hardy Space]
The {\it Hardy space} $H^2(\C^+)$ is the class of analytic functions $f$ on the upper half complex plane $\C^+:=\{ x+\i y \,|\, x \in \R, y>0\}$ which satisfies the growth condition
$$ \left( \int_{-\infty}^{\infty} |f(x+\i y)|^2 dx \right)^{1/2} < \infty,$$ for all $y>0$.
\end{Def}

Important properties of Hardy functions are given by Titchmarsh's Theorem \cite{N1972}.

\begin{Th}[Titchmarsh's Theorem] Let $g:=u+\i v \in H^2(\C^+)$. Then the following two assertions are equivalent:
\begin{itemize}
\item[1.] The Hardy function $g$ has no negative-frequency components. That is, $$g(z)={1 \over 2 \pi} \int_0^{\infty} G(\omega) e^{\i \omega z} d\omega,$$ where $G(\omega):= \int_{\R} g(x) e^{-\i \omega x} dx$ is the Fourier transform of $g$.
\item[2.] The real and imaginary parts verify the formulas: $$u(x+\i y)= u \ast P_y(x)=\int_{\R} P_y(x-t) u(t) dt,$$ and $$v(x+\i y)=v \ast Q_y(x)=\int_{\R} Q_y(x-t)v(t) dt,$$ for all $y>0$, where $P_y(x)={1 \over \pi} {y \over x^2+y^2}$ and $Q_y(x)={1 \over \pi} {x \over x^2+y^2}$ are the Poisson and conjugate Poisson kernel in $\C^+$.
\end{itemize}
\end{Th}
In this way, an analytic signal $f_A =f+\i \mathcal{H}[f]$ represents the boundary values of Hardy function $u+\i v$ in the upper half plane $\C^{+}$ \cite{G1981}. That is, $$f(x)=\lim_{y \rightarrow 0} u(x+\i y)$$ and $$\mathcal{H}[f](x) =\lim_{y \rightarrow 0} v(x+\i y).$$
Starting from this concept we are going to study the higher dimensional generalization on Clifford algebra.

\subsection{Clifford Algebra}
For all what follows we will work in ${\it Cl}_{0, m}$, the real {\it Clifford algebra}. Most of the basic knowledge and notations in relation to
Clifford algebra are referred to \cite{BDS, DSS}.  Let ${\bf e}_1, ..., {\bf e}_m $ be {\it basic elements} satisfying
${\bf e}_i{\bf e}_j+{\bf e}_j{\bf e}_i=-2\delta_{ij}$, where
$\delta_{ij}=1$ if $i=j,$ and $\delta_{ij}=0$ otherwise, $i, j=1, 2,
\cdots, m.$ The Clifford algebra ${\it Cl}_{0, m}$ is the associative algebra over
the real field $\R$. A general element in ${\it Cl}_{0, m}$, therefore, is of the form $x=\sum_S x_S {\bf e}_S,
x_s\in {\R}$, where ${\bf e}_S={\bf e}_{i_1}{\bf
e}_{i_2}\cdots {\bf e}_{i_l},$ and $S$ runs over all the ordered
subsets of $\{1,2,\cdots,m\},$ namely
$S=\{1 \leq i_1 <i_2< \cdots < i_l \leq m \}, \quad 1\leq l \leq
m.$

Let \[\ {\bf R}^m =\{\underline{x} \; |\; \underline{x}=x_1 {\bf e}_1 + \cdots + x_m {\bf e}_m, x_j
\in \R, j=1, 2, \cdots, m \}\] be identical with the usual Euclidean
space and an element in ${\bf R}^m$ is called a {\it vector}. Moreover, let $$
 {\bf R}_1^m =\{x \; |\; x=x_0+\underline{x}, x_0 \in \R, \ux \in {\bf R}^m \}$$ be the {\it para-vector} space and an element in ${\bf R}_1^m $ is called a {\it para-vector}.
The multiplication of two para-vectors
$x_0+\ux=\sum_{j=0}^{m}x_j {\bf e}_j$ and $y_0+\uy=\sum_{j=0}^{m}y_j {\bf e}_j$
is given by
$(x_0+\ux)(y_0+\uy)=(x_0y_0+\ux\cdot \uy) +(x_0 \uy+y_0\ux)+(\ux\wedge\uy)$
with
$
\ux\cdot \uy=-<\ux, \uy>=-\sum_{j=1}^{m}x_jy_j
$
and
$\ux\wedge\uy=\sum_{i<j}{\bf e}_{ij}(x_iy_j-x_jy_i).$
Clearly, we have
\begin{eqnarray}\label{addeq1}
\ux\cdot \uy=\frac{1}{2}(\ux \uy+\uy\ux)
\end{eqnarray}
and
$$\ux\wedge\uy=\frac{1}{2}(\ux\uy-\uy\ux).$$

There are three parts of $(x_0+\ux)(y_0+\uy)$. We denote them as follows
\begin{itemize}
  \item the {\it scalar part}: $x_0y_0+\ux\cdot \uy={\rm{Sc}}[(x_0+\ux)(y_0+\uy)]$,
  \item the {\it vector part }: $x_0 \uy+y_0\ux={\rm{Vec}}[(x_0+\ux)(y_0+\uy)]$,
  \item the {\it bi-vector part }: $\ux\wedge\uy={\rm{Bi}}[(x_0+\ux)(y_0+\uy)]$.
\end{itemize}
In particular, we have $\ux^2=-<\ux, \ux>=-|\ux|^2=-\sum_{j=1}^{m}x_j^2, \mbox{ for }\ux\in \R^m$.

The conjugation and reversion of ${\bf e}_S $ are defined by
$\overline{\bf e}_{S}=\overline{\bf e}_{il}\cdots\overline{\bf
e}_{i1}$ and ${\overline{\bf e}_j=-{\bf e}_j}$, respectively. Therefore, the Clifford conjugate
of a para-vector $x_0+\ux$ is
 $\overline{x_0+\ux}=x_0-\ux$. It is
easy to verify that $0 \not= x_0+\ux \in {\bf R}_1^m$ implies
$$(x_0+\ux)^{-1} := \frac{\overline{x_0+\ux}}{|x_0+\ux|^2}.$$
The open ball with center $0$ and radius $r$ in ${\bf R}_1^m $
is denoted by $B(0, r)$ and the unit sphere in ${\bf R}_1^m $ is
denoted by $S^m$.

The natural inner product between $x$ and $y$ in ${\it Cl}_{0, m}$ is defined by $<x, y> :=\sum_Sx_S\overline{y_S},$  where $x :=\sum_Sx_S{\bf e}_S, x_S \in \R$
and $y :=\sum_Sy_S{\bf e}_S, y_S \in \R$. The norm associated with
this inner product is defined by $|x|=<x, x>^{1\over 2}=(\sum_S|x_S|^2)^{{1\over 2}}.$

Let $\Omega$ be an open subset of ${\bf R}_1^m$ with a piecewise smooth boundary. We say that function $f$ defined on $\Omega$ such that $f(x_0+\ux)=\sum_S f_S
(x_0+\ux) {\bf e}_S$ is a Clifford-valued function or, briefly, a ${\it Cl}_{0,m}$-valued function, where $f_S$ are real-valued functions defined in $\Omega$.

A possibility to generalize complex analytic is offered by following the Riemann approach, which is introduced by means of the \emph{generalized Cauchy-Riemann
operator} $\frac{\partial}{\partial {x_0}}+\underline{D}$, where
$\underline D={\partial\over \partial x_1}{\bf e}_1+\cdots
+{\partial\over
\partial x_m}{\bf e}_m$ is the \emph{Dirac operator}. Nullsolutions to this operator provide us with the class of the so-called {\it monogenic functions}.

\begin{Def} (Monogenic Function)
An ${\it Cl}_{0, m}$-valued function $f$ is called left (resp. right) monogenic in $\Omega$ if $\left( \frac{\partial}{\partial {x_0}}+\underline{D}\right) f=0$ (resp. $f \left(\frac{\partial}{\partial {x_0}}+\underline{D} \right) =0$) in $\Omega$.
\end{Def}
In the following, let
$$ E(x_0+\ux)={\overline{x_0+\ux} \over |x_0+\ux|^{m+1}}$$ be the {\it Cauchy kernel} defined in ${\bf R}_1^{m}\setminus
\{0\}$. It is easy to verify that $E(x_0+\ux)$ is a monogenic function in ${\bf R}_1^{m}\setminus
\{0\}$ \cite{BDS, DSS}.
\begin{Rem}
\begin{itemize}
\item
For a ${\it Cl}_{0, m}$-valued function defined on an open subset of ${\bf R}^m$, we apply the {\it Dirac operator} $\underline{D}$ for the monogenic function.

\item
Throughout the paper, and unless otherwise stated, we only use left ${\it Cl}_{0, m}$-valued monogenic functions that, for simplicity, we call monoginic. Nevertheless, all results accomplished to left ${\it Cl}_{0, m}$-valued monogenic functions can be easily adapted to right ${\it Cl}_{0, m}$-valued monogenic functions.
\end{itemize}
\end{Rem}

We further introduce the right linear Hilbert space of integrable and square integrable ${\it Cl}_{0, m}$-valued functions in $\Omega \subset \R^m$ that we denote by  $L^1{(\Omega, {\it Cl}_{0, m})}$ and $L^2{(\Omega, {\it Cl}_{0, m})}$, respectively.
If $f \in L^1{(\R^m, {\it Cl}_{0, m})}$, the {\it Fourier
transform} of $f$ is defined by
\begin{eqnarray}\label{FT}\hat{f}(\uxi)=\int_{\R^m}e^{-\i<\ux,\uxi>}f(\ux)d\ux,\end{eqnarray}
if in addition, $\hat{f} \in L^1{(\R^m, {\it Cl}_{0, m})}$, then function $f$ can be recovered by the {\it inverse Fourier transform}
$$f(\ux)=\frac{1}{(2\pi)^m}\int_{\R^m}e^{\i<\ux,\uxi>}\hat{f}(\uxi)d\uxi.$$

The well-known Plancherel Theorem for Fourier transform of $f$ and $g \in L^2(\R^m, {\it Cl}_{0, m})$ holds
\begin{eqnarray*}\label{PT}
\int_{\R^m} f(\ux)g(\ux) d\ux=\frac{1}{(2\pi)^m}\int_{\R^m} \hat{f}(\uxi)\overline{\hat{g}(\uxi)} d\uxi.\end{eqnarray*}

In a recent paper \cite{FS1}, the authors defined the notion of the monogenic signal. It is regarded as an extension of the notion of the analytic signal to multidimensional signals.

\subsection{Monogenic Signal and Monogenic Scale Space}

The monogenic signal was defined by an
original signal and its "isotropic Hilbert transform" in the higher dimensional
spaces (a combination of the Riesz transforms).

\begin{Def}[Monogenic Signal]\label{MS} For $f \in L^2(\R^m, {\it Cl}_{0,m})$, the monogenic signal $f_M \in L^2(\R^m, {\it Cl}_{0,m})$ is defined by
$$f_{M}(\ux):=f(\ux)+H[f](\ux),$$ where $H[f]$ is the {\it isotropic Hilbert transform} of $f$
defined by
\begin{eqnarray*}
H[f](\ux)&:=& p.v.
\frac{1}{\omega_{m}}\int_{\R^m}\frac{\overline{\ux-\ut}}{|\ux-\ut|^{m+1}}f(\ut)d\ut\\
&=&\lim_{\epsilon\rightarrow
0^{+}}\frac{1}{\omega_{m}}\int_{|\ux-\ut|>\epsilon}\frac{\overline{\ux-\ut}}{|\ux-\ut|^{m+1}}f(\ut)d\ut\\
&=&-\sum_{j=1}^{m}R_j(f)(\ux)\e_j.
\end{eqnarray*}
Furthermore, $$R_j(f)(\ux):=\lim_{\epsilon\rightarrow
0^{+}}\frac{1}{\omega_{m}}\int_{|\ux-\ut|>\epsilon}\frac{x_j-t_j}{|\ux-\ut|^{m+1}}f(\ut)d\ut,$$
is the jth-Reisz transform of $f$ and $\omega_m=\frac{2\pi^{\frac{m+1}{2}}}{\Gamma(\frac{m+1}{2})}$ is the area of the unit sphere $S^m$ in $\R^m_1$.
\end{Def}

\begin{Rem} 
If $f(\ux)$ is real-valued, then by Definition \ref{MS}, $H[f](\ux)$ is
vector-valued.
\end{Rem}

Let us now generalize the notion of Hardy space to multidimensional space.

\begin{Def}[Monogenic Scale Space] \label{MSS}
The monogenic scale space $M^2(\R_1^{m,+})$ is the class of monogenic functions $f^+(\ux, s)$ defined on half space $$
 {\bf R}_1^{m,+} =\{x \; |\; x=(\underline{x}, s), \ux \in {\bf R}^m, s >0 \},$$  which satisfies the growth condition
$$ \left( \int_{\R^m} |f^+(\ux, s)|^2 d \ux \right)^{1/2} < \infty,$$
for all scale $s>0$.
\end{Def}

Like in the complex case, a monogenic signal is the boundary value of the monogenic scale function in the half space $\R_1^{m,+}$ \cite{LMcQ}. Some basic properties of the Monogenic scale space $M^2(\R_1^{m,+})$ in the half space are summarized as follows. For the proof of Theorem \ref{MT} we refer the reader to \cite{LMcQ} and \cite{KQ}.

\begin{Th}\label{MT}
Suppose $f^+(\ux, s) := u(\ux,s) + v(\ux,s) \in M^2(\R_1^{m,+})$. Then the following two assertions are equivalent:
\begin{itemize}
\item[1.] The inverse Fourier transform of $f^+$ vanishes for $s<0$. That is, the ${\it Cl}_{0,m}$-valued function $f^+(\ux, s)$ has the form
\begin{eqnarray*}\label{eq1}
f^{+}(\ux, s)=\frac{1}{(2\pi)^m}\int_{\R^m}e^{+}(s+\ux,
\ut)\hat{f}(\ut)d\ut
\end{eqnarray*}
 in the half space $s >0$, where $$e^{+}(s+\ux, \ut)=e^{-s|\ut|}e^{\i<\ux,
\ut>}(1+\i\frac{\ut}{|\ut|}),$$ is monogenic in ${\R}_1^m$ and $\hat{f}$ is the Fourier transform of $f$ given by (\ref{FT}).

\item[2.] The functions $u$ and $v$ are constructed by the Poisson and the conjugate Poisson integrals, respectively.
That is, \begin{eqnarray}\label{PK} u(\ux, s)=u*P_s(\ux)=\frac{1}{\omega_{m}}\int_{\R^m}\frac{s}{|s+(\ux-\ut)|^{m+1}} u(\ut)d\ut\end{eqnarray}
and \begin{eqnarray}\label{CPK}
\underline{v}(\ux, s)=v*Q_s(\ux)=\frac{1}{\omega_{m}}\int_{\R^m}
\frac{\overline{\ux-\ut}}{|s+(\ux-\ut)|^{m+1}}v(\ut)d\ut,\end{eqnarray} where
$P_s(\ux) :=\frac{1}{\omega_{m}}\frac{s}{|s+\ux|^{m+1}}$ and
$Q_s(\ux) :=\frac{1}{\omega_{m}}\frac{\overline{\ux}}{|s+\ux|^{m+1}}$ are the Poisson and
the conjugate Poisson kernel in $\R_1^{m,+}$, respectively.
\end{itemize}
\end{Th}
\section{Local Attenuation and Local Phase Vector}\label{S3}

Note that it is possible to write the monogenic scale function $f \in M^2(\R^{m,+}_1)$ in polar coordinate. Let us review the local feature \cite{YQS} as follows.

\begin{Def}[Local Features Representation I]
Suppose $f :=u+\uv \in M^2(\R^{m,+}_1)$ has the polar form \begin{eqnarray}\label{polar}
f(\ux, s) =A(f)(\ux, s)e^{\frac{\underline{v}(\ux, s)}{|\underline{v}(\ux,
s)|}\theta(\ux, s)},\end{eqnarray} then
\begin{eqnarray}\label{LA}
A(f)(\ux, s) : =|f(\ux, s)| =\sqrt{u(\ux, s)^2+|\underline{v}(\ux, s)|^2}\end{eqnarray} is
is called the {\it local amplitude}.
\begin{eqnarray}\label{PA}
\theta(\ux, s) :=\arctan \left(\frac{|\underline{v}(\ux, s)|}{u(\ux, s)}\right)\end{eqnarray}
is called the {\it phase angle} that is between $0$ and $\pi$.
\begin{eqnarray}\label{lpv}
\underline{r}(\ux, s) :=\frac{\underline{v}(\ux,
s)}{|\underline{v}(\ux, s)|}\theta(\ux, s),
\end{eqnarray}
 is called the {\it local phase vector}. ${\rm{Sc}}\left[(\underline{D}\theta(\ux,
s))\frac{\underline{v}(\ux, s)}{|\underline{v}(\ux, s)|}\right]$ is called the {\it directional phase derivative} and $e^{\underline{r}(\ux, s)}$ is called the {\it phase
direction}.
The \emph{phase derivative} or {\it instantaneous frequency} is defined
by
\begin{eqnarray}\label{pd}
{\rm{Sc}}\left[\left({\underline{D}f(\ux, s)}\right)\left({f(\ux,
s)}\right)^{-1}\right].
\end{eqnarray}
\end{Def}

Building on the ideas of \cite{FS2}, we can have the alternative form.
\begin{Def}[Local Features Representation II]
For nontrivial function $f:=u+\uv \in M^2(\R^{m,+}_1)$, the local amplitude is nonzero. We can rewrite (\ref{polar}) as
\begin{eqnarray}\label{eq22}
f(\ux, s) =e^{a(\ux, s)+\underline{r}(\ux, s)},
\end{eqnarray} where
\begin{eqnarray}\label{atten}
a(\ux, s) :=\ln A(f)(\ux, s)=\frac{1}{2}\ln(u(\ux,
s)^2+|\underline{v}(\ux, s)|^2)
\end{eqnarray} is called the {\it local attenuation}.
 \end{Def}

\begin{Rem} In one-dimensional case,
$\frac{\underline{v}(x, s)}{|\underline{v}(x, s)|}=\i,$ therefore the local phase vector
$\underline{r}(x, s)=\i \theta(x, s)$.
\end{Rem}

Suppose $f(x,s):=u(x,s)+\i v(x,s) \in H^2({\C^+})$ has the form (\ref{eq22}). That is, $f(x,s)=e^{a(x,s)+\i \theta(x,s)}$ has no zeros and isolated singularities in the half plane $\C^+$, then the local attenuation $a(x,
s)=\frac{1}{2}\ln(u^2+v^2)$ and the local phase
$\theta(x,s)=\arctan \left( \frac{v}{u}\right)$ are related by the Cauchy-Riemann system. The reason is that the composition of analytic function is analytic. If $f(x, s)=u(x, s)+ \i v(x, s)$
is analytic and has no zeros and isolated singularities in the half plane $\C^+$, then $a(x, s)+\i \theta(x, s)$ is also analytic in $\C^+$.

Using the Cauchy-Riemann system for $a(x, s)+\i \theta(x, s)$, we
have
\begin{eqnarray*}\label{eq3}
\frac{\partial a}{\partial s}+\frac{\partial \theta}{\partial x}=0,
\end{eqnarray*}
\begin{eqnarray*}\label{eq4}
\frac{\partial a}{\partial x}-\frac{\partial \theta}{\partial s}=0
\end{eqnarray*}
From the above system, we notice that:
\begin{itemize}
\item The instantaneous frequency $\frac{\partial \theta}{\partial x}$ can be
obtained by the minus of the scale derivative of the local
attenuation ${\partial a \over \partial s}$.
\item The zero points of the scale derivative of the local
phase ${\partial \theta \over \partial s}$ is given by the extrema of the local attenuation.\end{itemize}

Building on the ideas of 1D, the authors \cite{FS2} considered the {\it intrinsically 1D monogenic
signals}.
\begin{Def} If $f(x, s):=u(x,s)+\i v(x,s) \in H^2(\C^{+})$ has no zeros and isolated singularities in the half plane $\C^+$, then the {\it intrinsically 1D monogenic signal} is defined by
\begin{eqnarray}\label{1DMS}
f(<\ux, \un>, s)
&=&u(<\ux, \un>, s)+\overline{\un}v(<\ux, \un>, s)\nonumber\\
&=&u(<\ux, \un>, s)+\underline{v}(<\ux, \un>,
s\nonumber\\&=&e^{a(<\ux, \un>, s)+\ur(<\ux, \un>, s)},
\end{eqnarray}
where $\ux, \un \in \R^{m}$ and $\un$ is a fixed unit vector. The local attenuation is given by $a(<\ux, \un>, s)=\frac{1}{2}\ln(u^2+v^2)$ and the local phase vector is given by
$\ur(<\ux, \un>, s)=\overline{n}\arctan \left(\frac{v}{u}\right)$ for intrinsically 1D signal.
\end{Def}

Felsberg et al. \cite{FS2} proved that for the intrinsically 1D signals, the local
attenuation $a(<\ux, \un>, s)$ and the local phase-vector $\ur(<\ux,
\un>, s)$ are related by the Hilbert transform pairs (\ref{HTP}). Moreover, by the analyticity, using the generalized Cauchy-Riemann operator ${\partial \over \partial s}+\underline{D}$ on $a(<\ux, \un>,
s)+\underline{r}(<\ux, \un>, s)$, we have
\begin{eqnarray}\label{1Da}
\frac{\partial a}{\partial s}+\underline{D}(\underline{r})=0,
\end{eqnarray}
\begin{eqnarray}\label{1Db}
\underline{D}a+\frac{\partial \underline{r}}{\partial s}=0.
\end{eqnarray}

In \cite{FS2}, the {\it local frequency } of the
intrinsically 1D signal and the
{\it differential phase congruency} are defined by $\underline{D}(\underline{r})$ and $\frac{\partial
\underline{r}}{\partial s}$, receptively. From system (\ref{1Da}) and (\ref{1Db}), we notice that:\begin{itemize}
\item The local frequency
in an intrinsically 1D signal $\underline{D}(\ur)$ can also be obtained by the minus of the
scale derivative of the local attenuation ${\partial a \over \partial s}$. \item The zero points of the differential
phase congruency ${\partial r \over \partial s}$ is given by the extrema of the local attenuation.\end{itemize}

\begin{Rem}
In the recent paper \cite{YQS}, the instantaneous frequency of $f :=u+\uv=e^{a+\ur}$ is given by (\ref{pd})
\begin{eqnarray}\label{phd}
&&{\rm{Sc}}\left[({\underline{D}f(\ux, s)})({f(\ux,
s)})^{-1}\right]\nonumber\\
&=&{\rm{Sc}}\left[\left(\underline{D}\frac{\underline{v}}{|\underline{v}|}\right)\sin\theta(\ux,
s)\cos\theta(\ux,
s)\right]+{\rm{Sc}}\left[\left(\underline{D}\theta(\ux,
s)\right)\frac{\underline{v}}{|\underline{v}|}\right].\label{1}
\end{eqnarray}
In particular, if $\frac{\underline{v}}{|\underline{v}|}$ is a constant, the first term in (\ref{1}) vanishes,
then the instantaneous frequency is ${\rm{Sc}}\left[\left(\underline{D}\theta(\ux,
s)\right)\frac{\underline{v}}{|\underline{v}|}\right]$. It coincides with the local frequency defined in \cite{FS2}. That is, when
 $\frac{\underline{v}}{|\underline{v}|}=\overline{\un}$ is
a constant, the local frequency $\underline{D}(\underline{r})$ is given by
$(\underline{D}\theta(<\ux, \un>, s)) \overline{\un}$.
\end{Rem}

\begin{Rem} \begin{itemize}
\item In Clifford analysis \cite{DSS}, we notice that if $f(x,
s)=u(x, s)+\i v(x, s) \in H^2(\C^+)$, then for fixed unit
vector $\un \in \R^{m}$, the function
$$f(<\ux, \un>,
s)=u(<\ux, \un>, s)+\overline{\un}v(<\ux, \un>, s),$$ is monogenic in
$\C^+$. It is called {\it monogenic plane wave}.

\item Clearly, if $f(x, s)=e^{a(x,s)+i \theta (x,s)} \in H^2(\C^+)$ has no zeros and isolated singularities in $\C^+$, then
$a(x, s)+\i \theta(x, s)$ is also analytic in $\C^+$. Consequently, the function $a(<\ux, \un>, s)+\overline{\un}\theta(<\ux, \un>, s)=a(<\ux, \un>, s)+\underline{r} (<\ux, \un>, s)$
is monogenic in $\R^{m,+}_1$.
\end{itemize}
\end{Rem}

\begin{Prob}\label{P1} What is the situation in the higher dimension if the signal is not intrinsically 1D signal?\end{Prob}

The solution was not considered in \cite{FS2} and \cite{FDF}.
While in higher dimension, the situation is more complicated. The theory does not hold in general. 
In fact, if $f(\ux, s)=u(\ux, s)+\underline{v}(\ux, s)=e^{a(\ux, s)+\underline{r}(\ux, s)}$ is monogenic
in the half space $\R_1^{m,+}$, $a(\ux, s)+\underline{r}(\ux, s)$ is not monogenic in
general. Therefore, the local attenuation $a$ and the local phase vector $r$ are
not related by the generalized Cauchy-Riemann system in higher dimensions. Let us now look at an example to illustrate the topic discussed above.

\begin{Exa} Let $f(\ux,
s)=\frac{s}{|s+\ux|^{m+1}}+\frac{\overline{\ux}}{|s+\ux|^{m+1}}=E(s+\ux)$
be the Cauchy kernel in $\R_1^m \setminus \{0\}$, which is monogenic in ${\bf
R}_1^{m}\setminus \{0\}$. Then, by straightforward computations, we have
$$a(\ux, s)+\underline{r}(\ux,
s)=-\frac{m}{2}\ln(s^2+|\ux|^2)+\frac{\overline{\ux}}{|\ux|}\arctan \left(\frac{|\ux|}{s}\right).$$
Then we apply the generalized Cauchy-Riemann operator ${\partial \over \partial s}+\underline{D}$ on it, we have
\begin{eqnarray*}
&&\left(\frac{\partial}{\partial s}+\underline{D}\right)\left[-\frac{m}{2}\ln(s^2+|\ux|^2)+\frac{\overline{\ux}}{|\ux|}
\arctan \left(\frac{|\ux|}{s}\right)\right]\\
&=&\frac{(1-m)(s+\ux)}{s^2+|\ux|^2}+\frac{m-1}{|\ux|}\arctan \left(\frac{|\ux|}{s}\right)\neq 0.
\end{eqnarray*}
Therefore, $a(\ux, s)+\underline{r}(\ux, s)$ is not monogenic.
\end{Exa}

Let us now describe the solution for Problem \ref{P1}, Theorem \ref{th1} gives the relationship
between the local phase vector $r$ and the local attenuation $a$ in higher
dimensional spaces.

\begin{Th}\label{th1}
Let $f(\ux, s)=u(\ux, s)+\underline{v}(\ux, s)=e^{a(\ux,
s)+\underline{r}(\ux, s)} \in M^2(\R_1^{m,+})$, where $a(\ux, s)$ and $\ur(\ux, s)$ are the local
attenuation and the local
phase-vector defined by (\ref{atten}) and (\ref{lpv}), respectively. If $f$ has no zeros and isolated singularities in the half space $\R_1^{m,+}$. Then we have
\begin{eqnarray}\label{eq7}
\frac{\partial a}{\partial s}+{\rm
Sc}[(\underline{D}e^{\underline{r}})e^{-\underline{r}}]=0,
\end{eqnarray}
\begin{eqnarray}\label{eq8}
\frac{\partial \underline{r}}{\partial s}+\underline{D}a -{\rm
Vec} \left[ \left(\underline{D}\frac{\underline{v}}{|\underline{v}|} \right)\frac{\underline{v}}{|\underline{v}|} \right]\sin^2
\theta+(\sin \theta \cos \theta - \theta){\partial \over \partial s} \left({\uv \over |\ux|} \right)=0.
\end{eqnarray}
\end{Th}
In particular, if $\frac{\underline{v}}{|\underline{v}|}$ is independent of $s$, that is ${\partial \over \partial s} \left({\uv \over |\ux|} \right)=0$, then we
have the following corollary.

\begin{Cor}\label{cor1}
Let $f(\ux, s)=u(\ux, s)+\underline{v}(\ux, s)=e^{a(\ux,
s)+\underline{r}(\ux, s)}\in M^2(\R_1^{m,+})$, where $a$ and $\ur$ are the local
attenuation and local phase-vector defined by (\ref{atten}) and (\ref{lpv}), respectively.
If $f$ has no zeros and isolated singularities in the half space $\R_1^{m,+}$
and the local orientation $\frac{\uv}{|\uv|}$ does not change
through scale $s$, then we have
\begin{eqnarray}\label{eq88}
\frac{\partial \underline{r}}{\partial s}+\underline{D}a -{\rm
Vec}[(\underline{D}\frac{\underline{v}}{|\underline{v}|})\frac{\underline{v}}{|\underline{v}|}]\sin^2
\theta=0.
\end{eqnarray}
\end{Cor}

Combining (\ref{eq7}), (\ref{eq8}) and (\ref{phd}), we conclude that

\begin{Th}\label{th2} [Instantaneous Frequency]
\begin{itemize}
\item
The instantaneous frequency in higher dimensional spaces defined by
(\ref{pd}) is equal to the minus of the scale derivative of the
local attenuation ${\partial a \over \partial s}$.

\item The zero points of the differential phase congruency ${\partial r \over \partial s}$ is {\bf not} equal to the extrema of the local attenuation.
\end{itemize}
\end{Th}

\begin{Rem} By Theorem \ref{th2}, we notice that, like in one dimensional case,
 the phase derivative in higher dimensions
can also be given by the minus of the scale derivative of the local
attenuation. However, the zero points of
the phase congruency is {\bf not} equal to the extrema of the local
attenuation in high dimensional case. The nonzero extra term $$-{\rm
Vec} \left[ \left(\underline{D}\frac{\underline{v}}{|\underline{v}|} \right)\frac{\underline{v}}{|\underline{v}|} \right]\sin^2
\theta+(\sin \theta \cos \theta - \theta){\partial \over \partial s} \left({\uv \over |\ux|} \right)$$ appears in high dimensional cases.
\end{Rem}

We have divided the proof of Theorem \ref{th1} into a series of lemmas.

\begin{Lem}\label{lem2}
Let $f(\ux, s)=u(\ux, s)+\underline{v}(\ux, s)=e^{a(\ux,
s)+\underline{r}(\ux, s)} \in M^2(\R_1^{m,+})$, where $a(\ux, s)$ and $\ur(\ux, s)$ are the local
attenuation and the local
phase-vector defined by (\ref{atten}) and (\ref{lpv}), respectively. If $f$ has no zeros and isolated singularities in the half space $\R_1^{m,+}$. Then we have
\begin{eqnarray}\label{eq10}
{\rm{Sc}}\left[(\frac{\partial}{\partial s}e^{\underline{r}(\ux,
s)})e^{-\underline{r}(\ux, s)}\right]=0.
\end{eqnarray}
\end{Lem}

\noindent {\bf Proof: } By the generalized Euler formula
$e^{\underline{r}(\ux, s)}=e^{\frac{\underline{v}}{|\underline{v}|}\theta}=\cos\theta+\frac{\uv}{|\uv|} \sin
\theta$, we have
\begin{eqnarray}\label{qq}
&& \left(\frac{\partial}{\partial s}e^{\ur(\ux, s)}\right)  e^{-\ur(\ux, s)}\nonumber\\
&=&\frac{\partial}{\partial s} \left(\cos \theta+\frac{\uv}{|\uv|}\sin\theta\right) \left(\cos\theta-\frac{\uv}{|\uv|}\sin\theta\right) \nonumber\\
&=&\left(-\sin\theta \frac{\partial \theta}{\partial
s}+\frac{\partial \frac{\uv}{|\uv|}}{\partial
s}\sin\theta+\frac{\uv}{|\uv|}\cos\theta \frac{\partial
\theta}{\partial s}\right)
\left(\cos\theta-\frac{\uv}{|\uv|}\sin\theta\right)\nonumber\\
&=&\frac{\uv}{|\uv|} \frac{\partial\theta}{\partial
s}+\sin\theta\cos\theta \frac{\partial \frac{\uv}{|\uv|}}{\partial
s} -\sin^2\theta\frac{\partial \frac{\uv}{|\uv|}}{\partial
s}\frac{\uv}{|\uv|}.
\end{eqnarray}
Clearly, the scalar part of $\left(\frac{\partial}{\partial s}
e^{\ur(\ux, s)}\right) e^{-\ur(\ux, s)}$ is decided by the third
part of equation (\ref{qq}). Let us now prove the following
$$
{\rm{Sc}}\left[(\frac{\partial \frac{\underline{v}(\ux,
s)}{|\underline{v}(\ux, s)|}}{\partial s})\frac{\underline{v}(\ux,
s)}{|\underline{v}(\ux, s)|}\right] =0.
$$
Denote $\underline{I}(\ux, s):=\frac{\underline{v}(\ux, s)}{|\underline{v}(\ux, s)|}$, we have $\underline{I}^2(\ux, s)=-1$. Then $\frac{\partial [{\underline{I}^2(\ux, s)}]}{\partial s}=0$.
By equation (\ref{addeq1}), we have
\begin{eqnarray*}\label{addeq3}
\frac{\partial [{\underline{I}^2(\ux, s)}]}{\partial s}&=&\frac{\partial [{\underline{I}(\ux, s)}]}{\partial s}\underline{I}(\ux, s)+\underline{I}(\ux, s)\frac{\partial [{\underline{I}(\ux, s)}]}{\partial s}\nonumber\\
&=&2{\rm {Sc}}[\frac{\partial {\underline{I}(\ux, s)}}{\partial s}\underline{I}(\ux, s)]
=0.
\end{eqnarray*}
Thus, we obtain the desired result.

\begin{Lem}\label{lem3}
Let $f(\ux, s)=u(\ux, s)+\underline{v}(\ux, s)=e^{a(\ux,
s)+\underline{r}(\ux, s)} \in M^2(\R_1^{m,+})$, where $a(\ux, s)$ and $\ur(\ux, s)$ are the local
attenuation and the local
phase-vector defined by (\ref{atten}) and (\ref{lpv}), respectively. If $f$ has no zeros and isolated singularities in the half space $\R_1^{m,+}$. Then we have
\begin{eqnarray}\label{eq23}
{\rm{Vec}}\left[(\frac{\partial}{\partial s}e^{\underline{r}(\ux,
s)})e^{-\underline{r}(\ux,
s)}\right]=(\sin\theta\cos\theta-\theta)\frac{\partial\frac{\uv}{|\uv|}}{\partial
s}+\frac{\partial \ur}{\partial s}.
\end{eqnarray}
\begin{eqnarray}\label{eq24}
{\rm{Vec}}\left[(\underline{D}e^{\underline{r}(\ux,
s)})e^{-\underline{r}(\ux,
s)}\right]=-\sin^2\theta{\rm{Vec}}\left[(\underline{D}\frac{\uv}{|\uv|})\frac{\uv}{|\uv|}\right].
\end{eqnarray}
\end{Lem}

\noindent {\bf Proof: } From (\ref{qq}), we know that the vector part of
$(\frac{\partial}{\partial s}e^{\underline{r}(\ux,
s)})e^{-\underline{r}(\ux, s)}$ is decided by $\frac{\uv}{|\uv|}
\frac{\partial\theta}{\partial s}+\sin\theta\cos\theta
\frac{\partial \frac{\uv}{|\uv|}}{\partial s}$. Since $\underline{r}=\frac{\uv}{|\uv|}\theta$, we have
$$\frac{\partial \ur}{\partial s}=\frac{\partial\theta}{\partial
s}\frac{\uv}{|\uv|}+\theta\frac{\partial \frac{\uv}{|\uv|}}{\partial
s}.$$ Therefore, we obtain equation (\ref{eq23}).

To prove equation (\ref{eq24}), by direct calculation, we have
\begin{eqnarray}\label{qqq}
&& \left(\underline{D}e^{\ur(\ux, s)}\right)  e^{-\ur(\ux, s)}\nonumber\\
&=&\underline{D} \left(\cos\theta+\frac{\uv}{|\uv|}\sin\theta\right) \left(\cos\theta-\frac{\ur}{|\ur|}\sin\theta\right) \nonumber\\
&=&\left[-\sin\theta (\underline{D}\theta)+
(\underline{D}\frac{\uv}{|\uv|})\sin\theta+\cos\theta
(\underline{D}\theta)\frac{\uv}{|\uv|}\right]
\left[\cos\theta-\frac{\uv}{|\uv|}\sin\theta\right]\nonumber\\
&=&\frac{\uv}{|\uv|}(\underline{D}\theta)+\sin\theta\cos\theta
(\underline{D}\frac{\uv}{|\uv|})
-\sin^2\theta(\underline{D}\frac{\uv}{|\uv|})\frac{\uv}{|\uv|}.
\end{eqnarray}
The fist part and the second part of equation (\ref{qqq}) are scalar and bi-vector, respectively. Therefore the vector part of $(\underline{D}e^{\underline{r}(\ux,
s)})e^{-\underline{r}(\ux, s)}$ is decided by the third part of
equation (\ref{qqq}). Thus we obtain (\ref{eq24}).

We can now prove Theorem \ref{th1}.

\noindent {\bf Proof of Theorem \ref{th1}: } Since $f(\ux,
s)=e^{a(\ux, s)+\underline{r}(\ux, s)} \in M^2(\R_1^{m,+})$, we have
$$\left(\frac{\partial}{\partial s}+\underline{D}\right) e^{a(\ux,
s)+\underline{r}(\ux, s)}=0.$$ By straightforward computation, we have
$$e^{a(\ux, s)}\frac{\partial a(\ux, s)}{\partial s}e^{\underline{r}(\ux, s)}+e^{a(\ux, s)}\frac{\partial e^{\ur(\ux, s)}}{\partial s}
+e^{a(\ux, s)}[\underline{D}a(\ux, s)]e^{\ur(\ux, s)}+e^{a(\ux,
s)}(\underline{D}e^{\ur(\ux, s)})=0.$$ That is
\begin{eqnarray}\label{eq13}
\frac{\partial a(\ux, s)}{\partial s}+\frac{\partial e^{\ur(\ux,
s)}}{\partial s}e^{-\ur(\ux, s)} +\underline{D}a(\ux,
s)+(\underline{D}e^{\ur(\ux, s)})e^{-\ur(\ux, s)}=0.
\end{eqnarray}
Therefore, the scalar part of (\ref{eq13}) is zero. By combining Lemma \ref{lem2}, we have
\begin{eqnarray}\label{eq14}
&&{\rm Sc}\left[\frac{\partial a(\ux, s)}{\partial s}+\frac{\partial
e^{\ur(\ux, s)}}{\partial s}e^{-\ur(\ux, s)} +\underline{D}a(\ux,
s)+(\underline{D}e^{\ur(\ux, s)})e^{-\ur(\ux, s)}\right]\nonumber\\
&=&\frac{\partial a(\ux, s)}{\partial s}+{\rm Sc}[\frac{\partial
e^{\ur(\ux, s)}}{\partial s}e^{-\ur(\ux, s)}]+{\rm
Sc}[(\underline{D}e^{\ur(\ux, s)})e^{-\ur(\ux, s)}]\\
&=&\frac{\partial a(\ux, s)}{\partial s}+{\rm
Sc}[(\underline{D}e^{\ur(\ux, s)})e^{-\ur(\ux, s)}]
=0\nonumber.
\end{eqnarray}
Therefore, we get the desired result (\ref{eq7}).

The vector part of  Eq. (\ref{eq13}) is also zero. By using Lemma \ref{lem3}, we obtain
\begin{eqnarray}\label{eq15}
&&{\rm Vec}\left[\frac{\partial a(\ux, s)}{\partial
s}+\frac{\partial e^{\ur(\ux, s)}}{\partial s}e^{-\ur(\ux, s)}
+\underline{D}a(\ux, s)+(\underline{D}e^{\ur(\ux, s)})e^{-\ur(\ux,
s)}\right]\nonumber\\
&=&{\rm Vec}\left[\frac{\partial e^{\ur(\ux, s)}}{\partial
s}e^{-\ur(\ux, s)}\right] +\underline{D}a(\ux, s)+{\rm
Vec}\left[(\underline{D}e^{\ur(\ux, s)})e^{-\ur(\ux,
s)}\right]\\
&=&\frac{\partial \ur}{\partial
s}+\underline{D}a+(\sin\theta\cos\theta-\theta)\frac{\partial\frac{\uv}{|\uv|}}{\partial
s}-{\rm{Vec}}\left[(\underline{D}\frac{\uv}{|\uv|})\frac{\uv}{|\uv|}\right]\sin^2\theta
=0\nonumber.
\end{eqnarray}
This completes the proof.

If $f \in M^2(\R_1^{m,+})$ has the {\it axial form}
$$f(\ux, s)=u(\rho, s)+\frac{\overline{\ux}}{|\ux|}v(\rho, s),
\mbox{ }\rho=|\ux|.
$$ Then, in this case, the local orientation $\frac{\uv}{|\uv|}=\frac{\overline{\ux}}{|\ux|}$ does not change
through the scale $s$. By polar coordinate, $f(\ux, s)=e^{a(\rho,
s)+\frac{\overline{\ux}}{|\ux|}\theta(\rho, s)}$, using Theorem
\ref{th1}, we have the following corollary.

\begin{Cor}
Let $f(\ux, s)=u(\rho, s)+\frac{\overline{\ux}}{|\ux|}v(\rho,
s)=e^{a(\rho, s)+\frac{\overline{\ux}}{|\ux|}\theta(\rho, s)} \in M^2(\R_1^{m,+})$. Then we have
\begin{eqnarray}\label{eq11}
-\frac{\partial a}{\partial s}=\frac{\partial \theta}{\partial
\rho}+\frac{m-1}{\rho}\sin\theta\cos\theta,
\end{eqnarray}
\begin{eqnarray}\label{eq12}
\frac{\partial \theta}{\partial s}=\frac{\partial a}{\partial \rho}
+\frac{m-1}{\rho}\sin^2\theta.
\end{eqnarray}
\end{Cor}
It is easy to see that when $m=1$, the above system ((\ref{eq11}) and (\ref{eq12})) is just the Cauchy-Riemann system in the one dimensional case.

\section{Edge Detection Methods}\label{S4}

Edge detection by means of quadrature filters has two ways:
either by detecting local maxima of the local amplitude or by detecting
points of stationary phase in scale-space. In this section, we begin by reviewing the differential phase congruency method \cite{FS2}.

\subsection{Differential Phase Congruency Methods}
\begin{Method}[DPC] For intrinsically 1D monogenic signal $f \in H^2(\C^+)$ given by (\ref{1DMS}), if $f$ has no zero and isolated singularities in the half plane $\C^+$. Then the differential phase congruency (DPC) has the following formula
\begin{eqnarray}\label{ma}
\frac{\partial \underline{r}_{in1}(\ux, s)}{\partial s}=\frac{u(\ux,
s)\frac{\partial \underline{v}(\ux, s)}{\partial
s}-\underline{v}(\ux, s)\frac{\partial u(\ux, s)}{\partial s}
}{u(\ux, s)^2+|\underline{v}(\ux, s)|^2}=0,
\end{eqnarray} where $\underline{r}_{in1}(\ux, s):=\underline{r}(<\ux, \un>, s)$.
\end{Method}

By (\ref{1Db}), we notice that formula (\ref{ma}) can also be obtained by the $-\underline{D}a$. However, the zero points of the differential phase congruency is {\bf not} given by the extrema of the local attenuation in higher dimension.

\subsection{Proposed Methods}
Let's introduce the local attenuation (LA) method for monogenic signals.

\begin{Method}[LA]
For $f \in M^2(\R_1^{m,+})$ has no zeros and isolated singularities in the half space $\R^{m,+}_1$, the local attenuation has the formula
\begin{eqnarray}\label{LAM}
\underline{D}a(\ux, s)=\frac{u(\ux, s)\underline{D}[u(\ux,
s)]+|\underline{v}(\ux, s)|\underline{D}[|\underline{v}(\ux,
s)|]}{u^2(\ux, s)+|\underline{v}(\ux, s)|^2}.
\end{eqnarray}
\end{Method}

Applying Dirac operator $\underline{D}$ on the local attenuation $a$, by direct computation on (\ref{LA}), formula (\ref{LAM}) follows. Using (\ref{1Db}), we know that for intrinsically 1D signals,
the zero points of the differential phase congruency is given by the extrema of the local attenuation. Notice that formula (\ref{LAM}) is equivalent to (\ref{ma}) for the intrinsically 1D monogenic signal.

Our second method is the modified differential phase congruency (MDPC) method. To proceed, we need the following technical lemma.

\begin{Lem}\label{newth2}
\begin{eqnarray}\label{addeq5}
\frac{\partial \underline{r}(\ux, s)}{\partial
s}=\left(\theta-\sin\theta\cos\theta\right)\frac{\partial
\frac{\underline{v}(\ux, s)}{|\underline{v}(\ux, s)|}}{\partial s}+
\frac{u(\ux, s)\frac{\partial \underline{v}(\ux, s)}{\partial
s}-\underline{v}(\ux, s)\frac{\partial u(\ux, s)}{\partial s}
}{u^2(\ux, s)+|\underline{v}(\ux, s)|^2}.
\end{eqnarray}
\end{Lem}

\noindent{\bf Proof: } By  Eq. (\ref{lpv}), we have
\begin{eqnarray}\label{eq99}
\frac{\partial \underline{r}(\ux, s)}{\partial s}&=&\frac{\partial}
{\partial s}(\frac{\uv}{|\uv|}\theta)=\frac{\partial
\frac{\uv}{|\uv|}} {\partial
s}\theta+\frac{\uv}{|\uv|}\frac{\partial} {\partial s}\theta.
\end{eqnarray}
By straightforward computation, we have
\begin{eqnarray}\label{eq999}
\frac{\partial} {\partial s}\theta=\frac{\partial} {\partial
s}(\arctan \left(\frac{|\uv|}{u}\right))=\frac{\frac{\partial |\uv|}{\partial
s}u-|\uv|\frac{\partial u}{\partial s}}{u^2+|\uv|^2}.
\end{eqnarray}
 Then,
\begin{eqnarray}\label{eq100}
\frac{\uv}{|\uv|}\frac{\partial} {\partial
s}(\arctan \left(\frac{|\uv|}{u}\right))=\frac{\frac{\uv}{|\uv|}\frac{\partial
|\uv|}{\partial s}u-\uv\frac{\partial u}{\partial s}}{u^2+|\uv|^2}.
\end{eqnarray}
Using the equation
\begin{eqnarray*}
\frac{\partial \uv}{\partial s}=\frac{\partial}{\partial
s}(\frac{\uv}{|\uv|}|\uv|)
=|\uv|\frac{\partial \frac{\uv}{|\uv|}}{\partial
s}+\frac{\uv}{|\uv|}\frac{\partial |\uv|}{\partial s},
\end{eqnarray*}
we obtain
\begin{eqnarray}\label{eq101}
\frac{\uv}{|\uv|}\frac{\partial |\uv|}{\partial s}=\frac{\partial
\uv}{\partial s}-|\uv|\frac{\partial \frac{\uv}{|\uv|}}{\partial s}.
\end{eqnarray}
Applying  Eq. (\ref{eq101}) to  Eq. (\ref{eq100}), we have
\begin{eqnarray}\label{eq102}
\frac{\uv}{|\uv|}\frac{\partial} {\partial
s}(\arctan \left(\frac{|\uv|}{u}\right))=\frac{u\frac{\partial \uv}{\partial
s}-\uv\frac{\partial u}{\partial
s}}{u^2+|\uv|^2}-\frac{u|\uv|}{u^2+|\uv|^2}\frac{\partial
\frac{\uv}{|\uv|}}{\partial s}.
\end{eqnarray}
Combining  Eq. (\ref{eq99}) and  Eq. (\ref{eq102}), we obtain  Eq. (\ref{addeq5}).

\begin{Rem} Note that equation (\ref{ma}) is a
special case of (\ref{addeq5}). The reason is in the
intrinsically 1D neighborhood, the local orientation
$\frac{\underline{v}(\ux, s)}{|\underline{v}(\ux, s)|}=\overline{\un}$ is a
constant. So $\frac{\partial \frac{\underline{v}(\ux,
s)}{|\underline{v}(\ux, s)|}}{\partial s}=0$. In fact, formula (\ref{ma}) always holds if the local orientation is independent of
$s$.
\end{Rem}

Let us now define the points of modified differential phase congruency.
\begin{Def} Let $\underline{r}(\ux, s)$ be the local phase vector, given by (\ref{lpv}),
of function $f \in M^2(\R_1^{m,+}$. Points where
$$\frac{\partial \underline{r}(\ux, s)}{\partial s}-{\rm
Vec}[(\underline{D}\frac{\underline{v}}{|\underline{v}|})\frac{\underline{v}}{|\underline{v}|}]\sin^2
\theta+(\sin \theta \cos \theta - \theta)\frac{\partial
\frac{\underline{v}}{|\underline{v}|}}{\partial s}=0$$ are called
points of modified differential phase congruency (MDPC).
\end{Def}

\begin{Rem} From Theorem \ref{th1} we know that in
any higher dimensional cases, edge detection by means of local
amplitude maxima is equivalent to edge detection by modified
differential phase congruency.
\end{Rem}

Using Eq. (\ref{addeq5}), we can now proposed our second method, the so-called modified differential phase congruency (MDPC) method.

\begin{Method}[MDPC]
For $f \in M^2(\R_1^{m,+})$ has no zeros and isolated singularities in the half space $\R^{m,+}_1$, the MDPC has the formula
\begin{eqnarray}\label{eq9999}
&&\frac{\partial \underline{r}(\ux, s)}{\partial s}-{\rm
Vec}[(\underline{D}\frac{\underline{v}}{|\underline{v}|})\frac{\underline{v}}{|\underline{v}|}]\sin^2
\theta+\left(\sin \theta \cos \theta - \theta\right)\frac{\partial
\frac{\underline{v}}{|\underline{v}|}}{\partial s}\nonumber\\
&=&\frac{u(\ux, s)\frac{\partial \underline{v}(\ux, s)}{\partial
s}-\underline{v}(\ux, s)\frac{\partial u(\ux, s)}{\partial s}
}{u^2(\ux, s)+|\underline{v}(\ux, s)|^2}-{\rm
Vec}\left[(\underline{D}\frac{\underline{v}}{|\underline{v}|})\frac{\underline{v}}{|\underline{v}|}\right]\sin^2
\theta.
\end{eqnarray}
\end{Method}

Finally, we introduce a mixed method by combining local attenuation and modified differential phase congruency
(LA+MDPC) for edge detection.

\begin{Method}[LA+MDPC]
For $f \in M^2(\R_1^{m,+})$ has no zeros and isolated singularities in the half space $\R^{m,+}_1$, the MDPC has the formula
\begin{eqnarray}\label{eq.mix}
&&\frac{\partial \underline{r}}{\partial s}-\underline{D}a -{\rm
Vec}[(\underline{D}\frac{\underline{v}}{|\underline{v}|})\frac{\underline{v}}{|\underline{v}|}]\sin^2
\theta+(\sin \theta \cos \theta - \theta)\frac{\partial
\frac{\underline{v}}{|\underline{v}|}}{\partial s}\nonumber\\
&=&\frac{u(\ux, s)\frac{\partial \underline{v}(\ux, s)}{\partial
s}-\underline{v}(\ux, s)\frac{\partial u(\ux, s)}{\partial s}
}{u^2(\ux, s)+|\underline{v}(\ux, s)|^2}-\underline{D}a -{\rm
Vec}[(\underline{D}\frac{\underline{v}}{|\underline{v}|})\frac{\underline{v}}{|\underline{v}|}]\sin^2
\theta
\end{eqnarray}
\end{Method}

\section{Experiments}\label{S5}
In this section, we begin by showing the details of our proposed methods. Two classical edge detection methods, such as Canny and Sobel edge detectors, will be compared with our algorithms. The Canny edge detector will begin by applying Gaussian filter to the test images. Then Canny edge detector computes the gradients on the images.
For the Sobel edge detector, we only apply its gradients to the original test images.
For the DPC and our proposed methods, we first apply the Poisson filter to the test images, then we compute apply their formulas to the images.

By comparing with the classical methods, phased based methods may show more detail on image.

\begin{figure}[!]
  \centering
 \includegraphics[height=3.2cm]{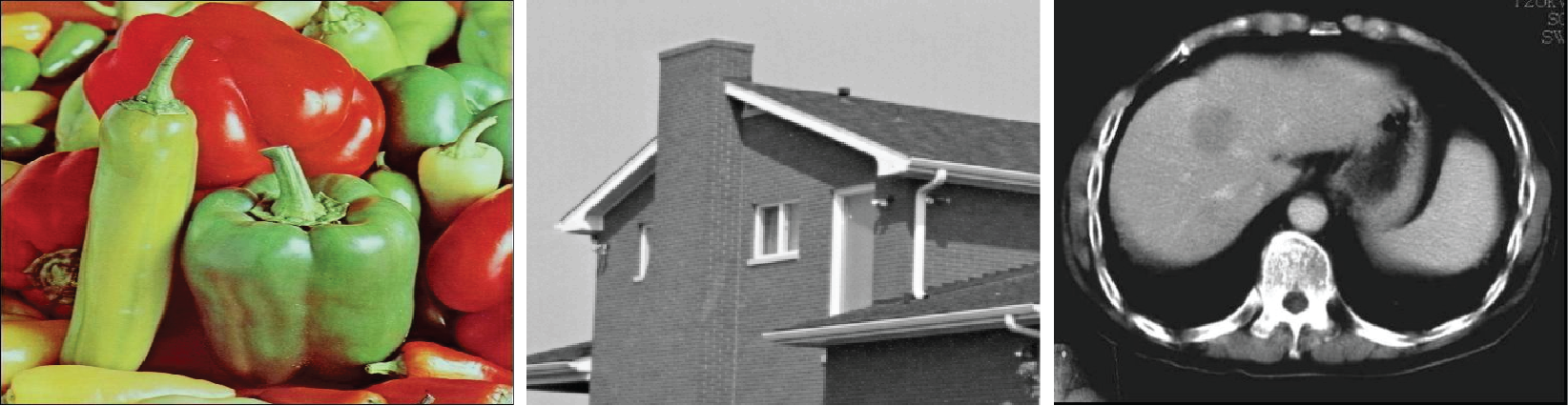}
  \caption{Original images}
  \label{fig.Original}
\end{figure}

\subsection{Algorithms}\label{S5.1}
Let us now give the details of the phase based algorithms. They are divided by the following steps.

\begin{description}
  \item[Step 1.]
  Input image $f(\ux)$.  For simplicity, the color image is converted to the gray image.

  \item[Step 2.]
   Poisson filtering: $u(\ux, s)=f*P_s(\ux)$ and and $\underline{v}(\ux, s)=f*Q_s(\ux)$  for a fixed scale $s>0$.
   We will discuss how to choose $s$ in Section \ref{5.1}. We consider $s$ in $0.1$, $0.5$,  $1.0$ and $5.0$. Moreover, we choose $s=0.5$ for all test images to compare with different methods.

  \item[Step 3.]
  Compute the local attenuation $a(\ux, s)=\frac{1}{2}\ln(u(\ux,s)^2+|\underline{v}(\ux, s)|^2)$ and
 local phase vector $\underline{r}(\ux, s)=\frac{\underline{v}(\ux,s)}{|\underline{v}(\ux, s)|}\theta(\ux, s)$, where
 the phase angle is given by $\theta(\ux, s)=\arctan \left(\frac{|\underline{v}(\ux, s)|}{u(\ux, s)}\right)$.

  \item[Step 4.]
  Compute gradients by different methods to get the gradient maps.

\begin{itemize}

  \item
  The differential phase congruency (DPC) method:  compute $\frac{\partial \underline{r}_{in1}(\ux, s)}{\partial s}$ by the formula
$$
\frac{u(\ux,s)\frac{\partial \underline{v}(\ux, s)}{\partial
s}-\underline{v}(\ux, s)\frac{\partial u(\ux, s)}{\partial s}
}{u(\ux, s)^2+|\underline{v}(\ux, s)|^2}.
$$

  \item
The local amplitude (LA) method: compute  $\underline{D}a(\ux, s)$, where $\underline{D}$ is the sum for the derivatives of image in vertical and horizontal directions.
By theoretical analysis,  $\underline{D}a(\ux, s)$ can be computed by
$$\frac{u(\ux, s)\underline{D}[u(\ux,
s)]+|\underline{v}(\ux, s)|\underline{D}[|\underline{v}(\ux,
s)|]}{u^2(\ux, s)+|\underline{v}(\ux, s)|^2}.$$

  \item
  The modified differential phase congruency (MDPC) method:  compute
  $\frac{\partial \underline{r}(\ux, s)}{\partial s}-{\rm
Vec}[(\underline{D}\frac{\underline{v}}{|\underline{v}|})\frac{\underline{v}}{|\underline{v}|}]\sin^2
\theta+\left(\sin \theta \cos \theta - \theta\right)\frac{\partial
\frac{\underline{v}}{|\underline{v}|}}{\partial s}$,
which equals to
$$\frac{u(\ux, s)\frac{\partial \underline{v}(\ux, s)}{\partial
s}-\underline{v}(\ux, s)\frac{\partial u(\ux, s)}{\partial s}
}{u^2(\ux, s)+|\underline{v}(\ux, s)|^2}-{\rm
Vec}\left[(\underline{D}\frac{\underline{v}}{|\underline{v}|})\frac{\underline{v}}{|\underline{v}|}\right]\sin^2
\theta.$$

  \item
  The mixed method by using local attenuation and modified differential phase congruency (LA+MDCP):  compute
  $$\frac{u(\ux, s)\frac{\partial \underline{v}(\ux, s)}{\partial s }-\underline{v}(\ux, s)\frac{\partial u(\ux, s)}{\partial s}}{u^2(\ux, s)+|\underline{v}(\ux, s)|^2}
  -\underline{D}a -{\rm Vec}[(\underline{D}\frac{\underline{v}}{|\underline{v}|})\frac{\underline{v}}{|\underline{v}|}]\sin^2\theta.$$
\end{itemize}

  \item[Step 5.]
 Applying Non-maximum suppress to these gradient maps,  which is the same as for the Canny edge detector.
After non-maximum suppression, the edge will become thinner  \cite{PK99, PK03}.
For a fair comparison, all the six methods aforementioned
utilize the non-maximum suppression method with the same parameters.
Concretely, we choose the radius $r=1.5$ and the lower and upper threshold values are $1.0$ and $3.5$, respectively.

\end{description}

\begin{figure}[!]
  \centering
 \includegraphics[height=20cm]{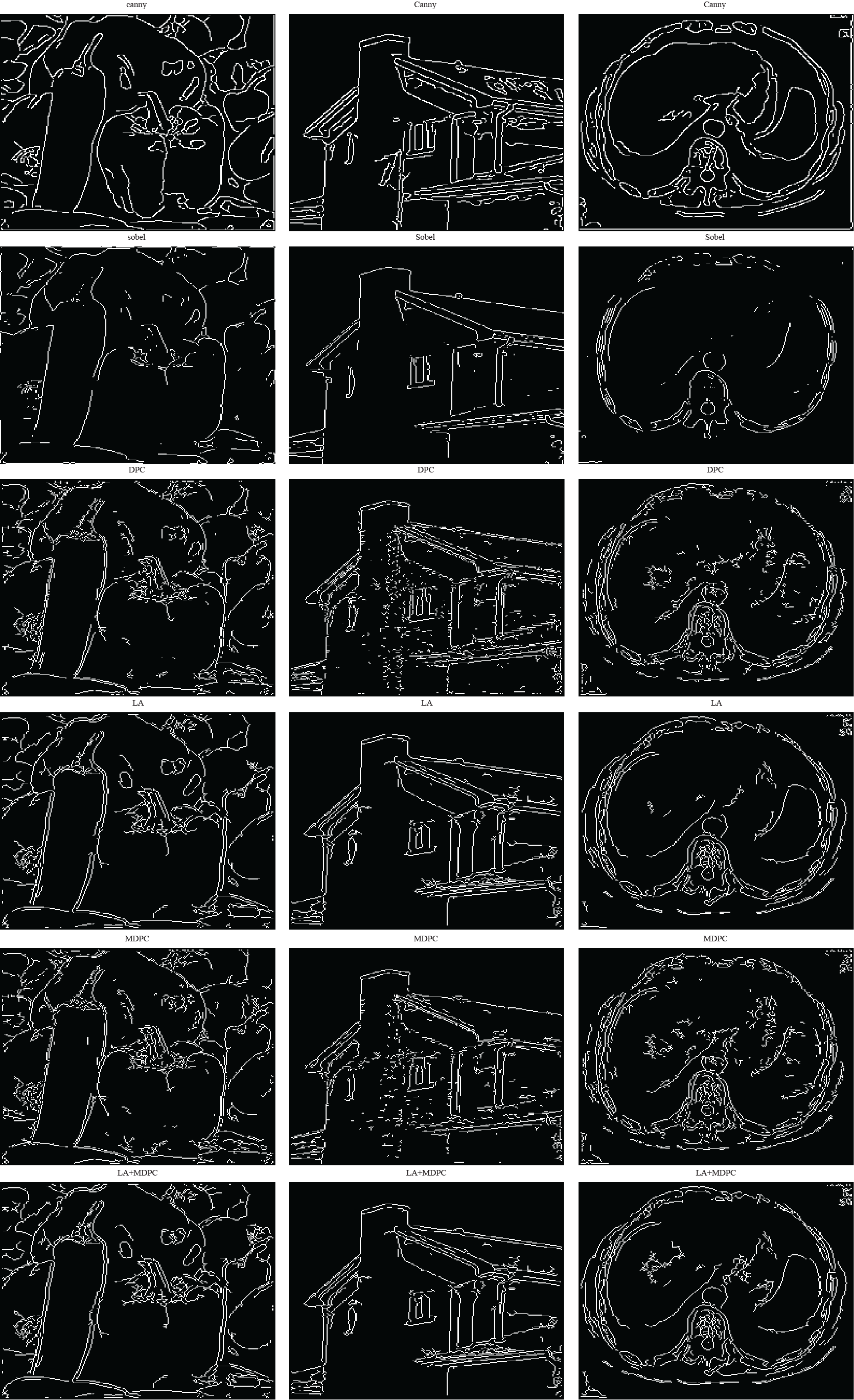}
  \caption{ Results for Canny, Sobel, DPC, LA, MDPC and LA+MDPC from top to the bottom.}
  \label{fig.3results}
\end{figure}

\subsection{Experiment Results}\label{S5.2}

We will use three different images (Fig. \ref{fig.Original}) for the comparison of different edge detectors.
Fig. \ref{fig.results} shows the edge detection results of various methods with the fixed scale $s=0.5$. From top to down of Fig. \ref{fig.3results}, there are six rows. Each row shows one comparison method. They are
Canny, Sobel, DPC, LA, MDPC and LA+MDPC methods, respectively. From the experiment results, we can draw the following conclusions.

\begin{itemize}
\item First, the mixed method yields decent edge detection results with fewer mistakes, outperforming other algorithms in some cases.

\item
Second, the comparison between the results of LA, DPC, MDPC and the mixed methods also
suggest that both local attenuation and local phase are important in edge detection.

\item
Our proposed method MDPC can achieve very good performances in dealing with the details. For the pepper in Fig. \ref{fig.3results}, we found that our method and canny's results are similar, we can find the edge of pepper in the results. However, for the shadows of the house and liver in Fig. \ref{fig.3results}, where the human eye is relatively subtle. Fortulately, the DPC and MDPC methods have found the details in the shadow. However,
 Canny, Sobel and LA methods cannot give the information about the shadows.
Canny uses the Gaussian filter which will make part of these shadows as noise and removed them.
While Sobel directly generate the horizontal and vertical differences, because the shadows
and the surrounding area is not much difference, which may not find the shadow of the image.
 By applying the phase based method, these details can be clearly found in our experiment results.
 This shows that our method can detect the whole smooth region and local small change region.
 Applications can be useful in places where it is difficult for the human eye to find the details.
\end{itemize}

\begin{Rem} For intrinsically 1D signals, we know
that edge detection by means of local amplitude maxima is equivalent
to edge detection by phase congruency. While, in intrinsically 2D
signals, we know that it dose not hold. In \cite{FS2}, the authors
said: $\lq\lq$ We cannot give an exhaustive answer to this question.
In this paper, since the behavior of phase and attenuation in
intrinsically 2D neighborhoods is still work in progress." From
 Eq. (\ref{eq9999}), we know that difference between the DPC and MDPC methods
 is ${\rm
Vec}\left[(\underline{D}\frac{\underline{v}}{|\underline{v}|})\frac{\underline{v}}{|\underline{v}|}\right]\sin^2
\theta$. By experiment, we know that the effect is not obvious.
\end{Rem}

\subsection{Effect of Scale}\label{S5.3}

\begin{figure}[!]
  \centering
 \includegraphics[height=13.5cm]{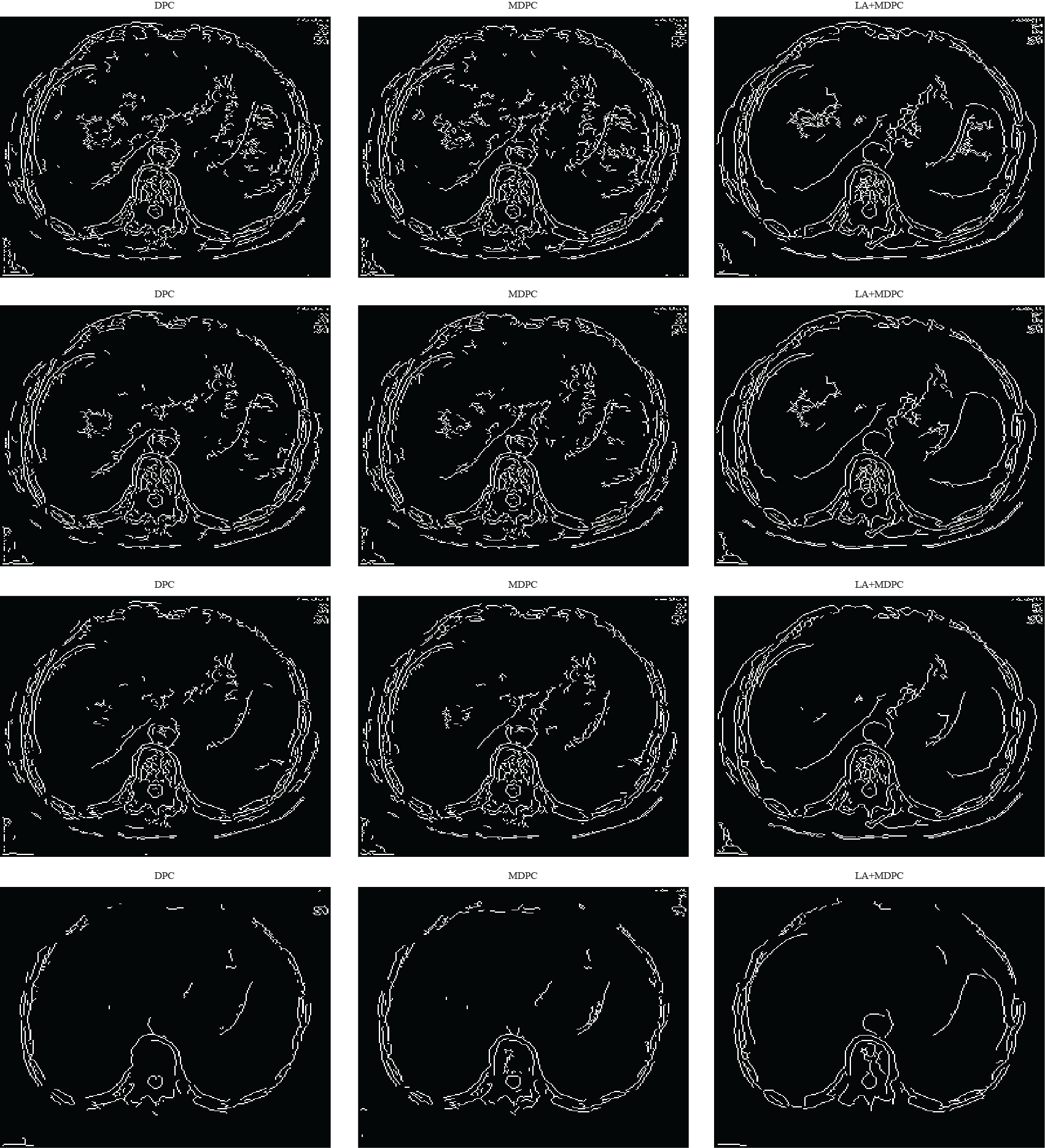}
  \caption{ Results for $s=0.1,~0.5,~1.0,~5.0$ from top to bottom.
 The first column show the differential phase congruency (DPC) method,
 the second column is the modified differential phase congruency (MDPC) method ,
 and the third column is the proposed mixed (LA+MDPC) method . }
  \label{fig.results}
\end{figure}

Monogenic signals at any scale $s>0$ form the monogenic scale-space $M^2(\R^{m,+}_1)$.
The representation of monogenic scale-space is just a monogenic function in the upper half space $\R^{m,+}_1$.
Therefore, considering the monogenic scale space instead of monogenic signal, it has a scale parameter $s>0$ to choose which provides us more analysis tools for different purposes.

We found that when $s$ tends to $0$, the Poisson integral tends to Hilbert transform.
Moreover, in the paper \cite{P}, Hilbert transform has been proved to be useful for image edge extraction.
In the choice of scale, we compare $s$ from $0.1$ to $5$, as can be seen in Fig. \ref{fig.3results},
when $s$ is smaller, the edge is more fine.
But too much detail is not good at all, so in the comparative experiment in Fig. \ref{fig.results},
we chose the case of 0.5 for $s$ to compare with the various methods.


\section{Acknowledgements}
 The authors acknowledge financial support from the National Natural Science Funds (No. 11401606) and University of Macau No. MYRG2015-00058-L2-FST and the Macao Science and Technology Development Fund FDCT/099/2012/A3.


\end{document}